# Bayesian inference with rescaled Gaussian process priors

**Aad van der Vaart and Harry van Zanten**[*]

*Department of Mathematics*
*De Boelelaan 1081a*
*1081 HV Amsterdam*
*e-mail:* aad@cs.vu.nl *and* harry@cs.vu.nl

**Abstract:** We use rescaled Gaussian processes as prior models for functional parameters in nonparametric statistical models. We show how the rate of contraction of the posterior distributions depends on the scaling factor. In particular, we exhibit rescaled Gaussian process priors yielding posteriors that contract around the true parameter at optimal convergence rates. To derive our results we establish bounds on small deviation probabilities for smooth stationary Gaussian processes.

## 1. Introduction

Gaussian processes have been adopted as building blocks for constructing prior distributions on infinite-dimensional statistical models in several settings. For instance, a sample path of a Gaussian process can be used directly as a prior model for a regression function (see e.g. [6], [22], [16]); after a monotonic transformation to the unit interval it can be used in the setting of classification (e.g. [16], [1], [5]); and after exponentiation and renormalization it becomes a model for density estimation (e.g. [11], [9], [18]).

Priors on functions of a single variable are commonly constructed using stationary Gaussian processes with smooth sample paths (e.g. [1], [5], [10], [16]). A popular example is the so-called *squared-exponential process*, i.e. the centered Gaussian process $W$ with covariance function $\mathbb{E} W_s W_t = a \exp(-b|t - s|^2)$ for some $a, b > 0$. The existing mathematical literature concerned with priors of this type focusses on computational issues (e.g. [16], [1], [10]) or posterior consistency ([5]). In the present paper we study posterior convergence rates, i.e. the rate at which the posterior distribution contracts around the true unknown functional parameter of interest. In particular, we are interested in exhibiting

[*]Partially supported by the Netherlands Organization for Scientific Research NWO

priors yielding optimal convergence rates if the unknown function belongs to a smoothness class.

It is well known that in the frequentist setup in which the data are sampled from a fixed true distribution, a prior on an infinite-dimensional model has to be carefully chosen in order to obtain optimal rates of contraction of the posterior (cf. e.g. [3], [17], [19], [4], [25], [23]). Even if posterior consistency is ensured, the rate of contraction of the posterior around the true functional parameter will typically be suboptimal if the regularity of prior process does not equal the regularity of the unknown parameter. Since Gaussian processes like the squared exponential process have infinitely often differentiable sample paths (at least in mean square sense), they will, at least without modification, typically not be appropriate as a prior model for a functional parameter with some finite smoothness level, in the sense that they yield suboptimal contraction rates.

In this paper we show however that this can be remedied by suitably rescaling the smooth process, with rescaling constants depending on the sample size. Given a fixed Gaussian process $(W_t : t \geq 0)$ indexed by the positive time axis and scaling constants $c_n > 0$ we use the rescaled sample path

$$t \mapsto W_{t/c_n}, \qquad t \in [0, 1] \tag{1.1}$$

as a prior model for a given function $w_0 : [0, 1] \to \mathbb{R}$ that indexes the distribution of the observations. The rescaling has the purpose of changing the appearance of the sample paths, so as to make them reflect more closely our prior ideas about the true parameter. Scaling factors $c_n \to \infty$ stretch the restrictions of the sample paths $t \mapsto W_t$ to the time interval $[0, 1/c_n]$ to the interval $[0, 1]$ and hence use only a small part of the randomness in the Gaussian process. Typically, this has the effect of smoothing the sample path. Conversely, scaling factors $c_n \to 0$ use the sample path $t \mapsto W_t$ on a long interval $[0, 1/c_n]$ and shrink this to the interval $[0, 1]$. This typically makes the prior rougher, by incorporating the randomness of a longer time period.

Coming back to the example of the squared exponential process, we will show that for any given regularity level $\alpha > 0$ there exist scaling factors $c_n \to 0$ ("roughening of the sample paths") such that after rescaling, we obtain a prior yielding (up to logarithmic factors) optimal contraction rates if the true parameter is $\alpha$-smooth. To prove this result we use the general theory on posterior convergence rates for Gaussian priors developed in Van der Vaart and Van Zanten [20]. The results of the latter paper state that the rate of convergence for a Gaussian process prior $W$ is determined by the behaviour of its *concentration function*

$$\varphi_{w_0}(\varepsilon) = \inf_{\|h - w_0\| \leq \varepsilon} \|h\|_{\mathbb{H}}^2 - \log \Pr\left(\|W\| \leq \varepsilon\right)$$

for $\varepsilon \to 0$, where $\mathbb{H}$ is the reproducing kernel Hilbert space (RKHS) associated to the process $W$, $\|\cdot\|_{\mathbb{H}}$ is the RKHS-norm, and $\|\cdot\|$ is the norm of the function space where $W$ takes its values. More precisely, the results state that asymptotically, the posterior concentrates its mass on balls around the true parameter with radius of the order $\varepsilon_n$, where $\varepsilon_n$ is found by solving

$$\varphi_{w_0}(\varepsilon_n) \lesssim n\varepsilon_n^2. \tag{1.2}$$

Our results for rescaled smooth, stationary Gaussian process priors are obtained by studying the RKHS and small deviations behaviour of such processes, leading to upper bounds for their concentration functions. For $W$ the centered Gaussian process with covariance function $\mathbb{E}W_s W_t = \exp(-|t-s|^2)$ and $c > 0$ we find that for the rescaled process $W^c = (W_{t/c} : 0 \le t \le 1)$, and $w_0$ a $C^\alpha$ function,

$$-\log \Pr\Big(\sup_{0 \le t \le 1} |W_t^c| \le 2\varepsilon\Big) \lesssim \frac{1}{c}\Big(\log \frac{1}{c\varepsilon^2}\Big)^2$$

and

$$\inf_{\|h-w_0\|_\infty \le c^\alpha} \|h\|_{\mathbb{H}^c}^2 \lesssim \frac{1}{c}.$$

This implies that if we use rescaling rates $c_n \to 0$, (1.2) is solved for $\varepsilon_n \gtrsim c_n^\alpha \vee (\log n)/\sqrt{nc_n}$. The optimal choice $c_n \sim (n/\log^2 n)^{-1/(1+2\alpha)}$ yields the rate $\varepsilon_n \sim (n/\log^2 n)^{-\alpha/(1+2\alpha)}$. Up to a logarithmic factor, this is the well-known minimax rate for estimating $\alpha$-regular functions.

In addition to smooth stationary prior processes we also consider self-similar processes. At first thought one might expect that rescaling would have no effect on such processes, but this turns out to be false. Stretching and shrinking causes the smoothing and roughening effects mentioned previously. Convergence rates for posteriors based on certain self-similar Gaussian process priors were obtained in [20]. We proved for instance that if $W$ is a $k$-fold integrated Brownian motion (plus an independent polynomial part), then using $W$ as prior on $(k+1/2)$-smooth functions yields an optimal convergence rate for the posterior. In this paper we show that after rescaling, this prior becomes appropriate for a larger range of smoothness levels. For any $\alpha \in (0, k+1]$ there exist scaling factors $c_n$ such that the prior based on the rescaled process $W_{t/c_n}$ yields optimal contraction rates if the true parameter is $\alpha$-smooth. In this case we have $c_n \to 0$ ("roughening") if $\alpha < k + 1/2$ and $c_n \to \infty$ ("smoothing") if $\alpha > k + 1/2$. The range of $\alpha$'s for which rescaling the $k$-fold integrated Brownian motion leads to a rate-optimal posterior is limited by the smoothness level $k+1$ of the functions in the RKHS of the process. Technically, the results for self-similar processes are relatively easy consequences of the general results obtained in [20].

The results of this paper can be viewed as mathematical support for the common use of rescaled Gaussian process priors in Bayesian practice (see for instance [1], [10], [24]). We show that, from a frequentist perspective, rescaling greatly enlarges the range of models for which a given Gaussian process prior is appropriate. In a practical setting one often tries to robustify a Bayes procedure, or reduce subjectivity, by employing a random rescaling, i.e. using the prior $W_{t/C}$ with a random scaling factor $C$ independent of $W$, rather than the prior $W_t$ itself. Further analysis is necessary on this issue, but the results in this paper may serve as a starting point for such an investigation.

Another point that deserves elaboration is the extension of our results to multivariate settings, i.e. cases where the unknown function of interest is a function of several variables. This requires however the generalization to higher dimensions of a number of approximation results (cf. Lemmas 2.2 and 2.3), which, in general, can be technically quite involved.

The remainder of this paper is organized as follows. In the next section we introduce and study the Gaussian processes that will serve as prior models. To prepare for the proofs of our main results on posterior convergence rates we obtain small deviation bounds and results on the approximation properties of the RKHS of rescaled smooth stationary processes and multiply integrated Brownian motions. The results on posterior contraction are precisely stated in the final Section 3.

The notation $\lesssim$ is used for "smaller than or equal to a universal constant times", and $\asymp$ is "proportional up to constants". We use the notation $C[0,1]$ for the space of continuous functions $f : [0,1] \to \mathbb{R}$, endowed with the uniform norm, and $L_r(\mu)$ for the measurable functions $f : [0,1] \to \mathbb{R}$ or $f : [0,1] \to \mathbb{R}$ with $\|f\|_r^r := \int |f|^r \, d\mu < \infty$. Furthermore, for $\beta > 0$ we let $C^\beta[0,1]$ denote the Hölder space of order $\beta$, consisting of the functions $f \in C[0,1]$ that have $\underline{\beta}$ continuous derivatives, for $\underline{\beta}$ the biggest integer strictly smaller than $\beta$, with the $\underline{\beta}$th derivative $f^{(\underline{\beta})}$ being Lipshitz continuous of order $\beta - \underline{\beta}$. For $\varepsilon > 0$ let $N(\varepsilon, B, d)$ be the minimum number of balls of radius $\varepsilon$ needed to cover a subset $B$ of a metric space with metric $d$.

## 2. Prior processes

The theorems on posterior contraction rates that we present in the next section concern two classes of priors. The first are constructed by rescaling smooth, stationary Gaussian processes, the second by rescaling multiply integrated Brownian motions. In this section we study these rescaled processes, obtaining results on their small deviations behaviour and the approximation properties of their reproducing kernel Hilbert spaces (RKHSs). Together with the general theory of [20], these results will allow us to obtain rates of convergence for posteriors.

### 2.1. Smooth stationary processes

Consider a centered, mean-square continuous Gaussian process $W = (W_t : t \geq 0)$ with covariance function

$$\mathbb{E}W_s W_t = \varphi(s-t), \tag{2.1}$$

for a given continuous function $\varphi : \mathbb{R} \to \mathbb{R}$. For a fixed scaling constant $c > 0$, we define the rescaled version $W^c$ of the process $W$ by setting $W_t^c = W_{t/c}$.

By Bochner's theorem the function $\varphi$ is representable as the characteristic function

$$\varphi(t) = \int e^{-it\lambda} \, d\mu(\lambda)$$

of a symmetric, finite measure $\mu$ on $\mathbb{R}$, called the *spectral measure* of the process $W$. (The minus sign in the exponent is for consistency in notation, but is superfluous as $\mu$ is symmetric.) We shall consider spectral measures satisfying the condition

$$\int e^{\delta|\lambda|} \, \mu(d\lambda) < \infty \tag{2.2}$$

for some $\delta > 0$. This condition on the tails of the spectral measure should be viewed as a smoothness condition on $W$. It implies for instance that the sample paths of $W$ are infinitely often differentiable, at least in mean-square sense. Examples of processes satisfying (2.2) are the centered Gaussian processes with covariance functions $(s, t) \mapsto \exp(-|t - s|^2)$ or $(s, t) \mapsto (1 + |t - s|^2)^{-1}$, which correspond, respectively, to Gaussian and Laplace spectral measures. Observe that in particular (2.2) implies that $\mu$ has finite moments of all orders and hence, since

$$\mathbb{E}(W_s - W_t)^2 = \int |e^{i\lambda t} - e^{i\lambda s}|^2 \, \mu(d\lambda) \leq |t - s|^2 \int \lambda^2 \, \mu(d\lambda),$$

the process $W$ admits a continuous version. We shall work with this continuous version throughout.

The spectral measure $\mu_c$ of the rescaled process $W^c$ is obtained by rescaling $\mu$:

$$\mu_c(B) = \mu(cB).$$

Let $L_2(\mu_c)$ be the set of all functions $h : \mathbb{R} \to \mathbb{C}$ whose modulus $|h|$ is square integrable with respect to $\mu_c$. Denote by $\mathcal{F}_c h$ the transform $\mathcal{F}_c h : \mathbb{R} \to \mathbb{C}$ of the function $h$ relative to the measure $\mu_c$:

$$(\mathcal{F}_c h)(t) = \int e^{-it\lambda} h(\lambda) \, d\mu_c(\lambda).$$

Note that $\mathcal{F}_c$ maps $L_2(\mu_c)$ into the space $C(\mathbb{R})$ of continuous functions on the real line.

The following lemma describes the RKHS $\mathbb{H}^c$ of the process $(W_t^c : t \in [0, 1])$. Recall that this space is defined as the completion of the linear span of the functions $h_s$ defined by $h_s(t) = \mathbb{E}W_s^c W_t^c$, with $s \in [0, 1]$, under the inner product $\langle h_s, h_t \rangle_{\mathbb{H}^c} = \mathbb{E}W_s^c W_t^c$.

**Lemma 2.1.** *Under condition (2.2) the reproducing kernel Hilbert space of the process $(W_t^c : 0 \leq t \leq 1)$ (viewed as a map in $C[0, 1]$) is the set of real parts of all transforms $\mathcal{F}_c h$ (restricted to the interval $[0, 1]$) of functions $h \in L_2(\mu_c)$, equipped with the square norm*

$$\|\mathcal{F}_c h\|_{\mathbb{H}^c}^2 = \|h\|_{L_2(\mu_c)}^2 = \int |h|^2 \, d\mu_c.$$

*Proof.* Although the RKHS is real by definition, it will be convenient to complete the linear span of the functions $h_s$ over the complex numbers. Because the functions $h_s$ are real-valued, the RKHS is the set of real parts of functions in this complex RKHS.

If $e_s : \mathbb{R} \to \mathbb{C}$ denotes the function $e_s(\lambda) = e^{is\lambda}$, then, by the definition of the spectral measure,

$$\langle h_s, h_t \rangle_{\mathbb{H}^c} = \mathbb{E}W_s^c W_t^c = \varphi\Big(\frac{s-t}{c}\Big) = \int e^{i(s-t)\lambda}\, d\mu_c(\lambda) = \langle e_s, e_t \rangle_{L_2(\mu_c)}.$$

It follows that the linear extension of the map $L : h_s \mapsto e_s$ is an isometry for the Hilbert space structures from $\operatorname{lin}(h_s : s \in [0,1]) \subset \mathbb{H}^c$ onto $\operatorname{lin}(e_s : s \in [0,1]) \subset L_2(\mu_c)$. Hence, $\mathbb{H}^c$ is the inverse image under $L$ of the closure of the set of functions $\operatorname{lin}(e_s : s \in [0,1])$ in $L_2(\mu_c)$. Now, again by the definition of the spectral measure,

$$\Big(L^{-1}e_s\Big)(t) = h_s(t) = \mathbb{E}W_s^c W_t^c = \varphi\Big(\frac{t-s}{c}\Big) = (\mathcal{F}_c e_s)(t).$$

It follows that the inverse $L^{-1}$ is exactly the transform $\mathcal{F}_c$.

Finally we prove that $\operatorname{lin}(e_s : s \in [0,1])$ is dense in $L_2(\mu_c)$, using the condition (2.2). As $s \downarrow 0$, by dominated convergence,

$$\frac{e_s - e_0}{s}(\lambda) \to i\lambda,$$

as functions in $L_2(\mu_c)$. Because the functions on the left side are in $\operatorname{lin}(e_s : s \in [0,1])$, the function $\lambda \mapsto i\lambda$ is in the closure of this set. We repeat this argument with the function $(e_s - e_0)(\lambda)/s - i\lambda$ to see that the function $\frac{1}{2}(i\lambda)^2$ is contained in the closure of $\operatorname{lin}(e_s : s \in [0,1])$, and so on. We conclude that all polynomials $\lambda \mapsto \lambda^k$ are contained in this closure. By extension to the complex case of Proposition 6.4.1 of [15] and (2.2), it then follows that this closure is dense in the full space $L_2(\mu_c)$. □

Let $\varphi_c(x) = \varphi(x/c)$ and denote by $\varphi_c * G(t) = \int \varphi_c(t-s)\, dG(s)$ the density of the convolution of a signed measure $G$ and the distribution corresponding to the density $\varphi_c$. By Fubini's theorem, such a convolution can be written as

$$\varphi_c * G(t) = 2\pi\, (\mathcal{F}_c \hat{G})(t),$$

for $\hat{G}$ the characteristic function of $G$, defined by $2\pi\hat{G}(\lambda) = \int e^{it\lambda}\, dG(t)$. Because $2\pi|\hat{G}|$ is uniformly bounded by the total variation of $G$, it is contained in $L_2(\mu_c)$ and hence the function $\varphi_c * G$ is contained in the RKHS, with square norm

$$\|\varphi_c * G\|_{\mathbb{H}^c}^2 = (2\pi)^2 \int |\hat{G}|^2\, d\mu_c.$$

For a measure $G$ on the interval $[0,1]$ this follows readily from the definition of the RKHS as the linear space spanned by the functions $t \mapsto \mathbb{E}W_s^c W_t^c = \varphi_c(s-t) = \varphi_c * \delta_s(t)$. As shown by the preceding lemma, under condition (2.2), the functions $\varphi_c * G$ are contained in the RKHS for any signed measure $G$ on the full line $\mathbb{R}$. This will be important for the proof of the following lemma, which quantifies how well $C^\beta$ functions can be approximated by elements of the RKHS of the rescaled process $W^c$.

**Lemma 2.2.** *Let $\mu$ satisfy (2.2) and possess a Lebesgue density that is bounded away from zero on a neighborhood of 0. Let $\beta > 0$ be given. Then for any $w \in C^\beta[0, 1]$ there exist constants $C_w$ and $D_w$ depending only on $w$ such that, as $c \downarrow 0$,*

$$\inf\big\{\|h\|_{\mathbb{H}^c}^2 : \|h - w\|_\infty \leq C_w c^\beta\big\} \leq D_w\big(\frac{1}{c}\big).$$

*Proof.* Let $\underline{\beta}$ be the biggest integer strictly smaller than $\beta$. Let $\hat{\varphi}$ be the density of $\mu$, and let $\hat{\psi}(\lambda) = (2\pi)^{-1} \int e^{it\lambda} \psi(t) \, dt$ be the Fourier transform of a general function $\psi : \mathbb{R} \to \mathbb{C}$. The Fourier transform of the function $\psi_c$ defined by $\psi_c(x) = \psi(x/c)$ is given by $\hat{\psi_c}(\lambda) = c\hat{\psi}(c\lambda)$.

There exists a symmetric, integrable function $\psi : \mathbb{R} \to \mathbb{R}$ with $\int \psi(t) \, dt = 1$, $\int t^k \psi(t) \, dt = 0$ for every $k = 1, \ldots, \underline{\beta}$, and $\int |t|^\beta |\psi|(t) \, dt < \infty$, and such that the function $|\hat{\psi}|^2/\hat{\varphi}$ is bounded. Take for instance a function $\psi$ with compactly supported, symmetric, real-valued Fourier transform $\hat{\psi}$ which equals $1/(2\pi)$ in a neighborhood of zero, so that

$$\frac{1}{2\pi} \int (it)^k \psi(t) \, dt = \hat{\psi}^{(k)}(0) = \begin{cases} 0, & k \geq 1, \\ \dfrac{1}{2\pi}, & k = 0. \end{cases}$$

We can extend $w : [0, 1] \to \mathbb{R}$ to a function $w : \mathbb{R} \to \mathbb{R}$ with compact support and $\|w\|_\beta < \infty$.

By Taylor's theorem we can write, for $s, t \in \mathbb{R}$,

$$w(t + s) = \sum_{j=0}^{\underline{\beta}} w^{(j)}(t) \frac{s^j}{j!} + S(t, s),$$

where, for some $\xi \in [0, 1]$,

$$|S(t, s)| = \frac{|s|^{\underline{\beta}}}{\underline{\beta}!} \big| w^{(\underline{\beta})}(t + \xi s) - w^{(\underline{\beta})}(t) \big| \leq \frac{|s|^\beta}{\underline{\beta}!} \|w\|_\beta.$$

In view of the assumption that $\psi$ is a higher order kernel, for any $t \in \mathbb{R}$,

$$\frac{1}{c}(\psi_c * w)(t) - w(t) = \int \psi(s)\big(w(t - sc) - w(t)\big) \, ds = \int \psi(s) S(t, -cs) \, ds.$$

Combining the preceding displays shows that $\|c^{-1}\psi_c * w - w\|_\infty \leq c^\beta \|w\|_\beta K/\underline{\beta}!$ for $K = \int |s|^\beta |\psi|(s) \, ds$.

For $\hat{w}$ the Fourier transform of $w$, we can write

$$(w * \psi_c)(t) = 2\pi \int e^{-it\lambda} \hat{w}(\lambda) \hat{\psi_c}(\lambda) \, d\lambda = 2\pi \mathcal{F}_c\big(\frac{\hat{w}\hat{\psi_c}}{\hat{\varphi_c}}\big)(t).$$

It follows that the function $c^{-1}\psi_c * w$ is contained in the RKHS, with square norm

$$\frac{1}{c^2}\|w * \psi_c\|_{\mathbb{H}^c}^2 = \frac{1}{c^2}(2\pi)^2 \int \left|\frac{\hat{w}\hat{\psi}_c}{\hat{\varphi}_c}\right|^2 d\mu_c = \frac{1}{c}(2\pi)^2 \int \frac{|\hat{w}(\lambda)|^2 |\hat{\psi}(c\lambda)|^2}{\hat{\varphi}(c\lambda)} d\lambda$$

$$\leq \frac{1}{c}(2\pi)^2 \int |\hat{w}(\lambda)|^2 d\lambda \left\|\frac{|\hat{\psi}|^2}{\hat{\varphi}}\right\|_\infty.$$

Here $(2\pi)^2 \int |\hat{w}(\lambda)|^2 d\lambda = \int w^2(t) dt$ is finite. □

Lemma 2.1 implies that under (2.2) the elements of the RKHS can be continuously extended to functions that are analytic on the strip $\{t \in \mathbb{C} : |\operatorname{Im} t| < (c\delta)/2\}$ in the complex plane. The following entropy estimate is therefore related to classical results on the entropy of spaces of analytic functions, obtained by Kolmorogov and Tihomirov [7]. They obtain estimates of the order $\left(\log(1/\varepsilon)\right)^2$ for the entropy of the spaces we are interested in. The following lemma makes the dependence on the scaling constant $c$ explicit, which is essential for the proof of our main results.

**Lemma 2.3.** *Assume that the spectral measure satisfies (2.2) for some $\delta > 0$. Then the entropy of the unit ball $\mathbb{H}_1^c$ of the RKHS of the process $W^c = (W_t^c : 0 \leq t \leq 1)$ (viewed as map in $C[0,1]$) satisfies*

$$\log N\big(\varepsilon, \mathbb{H}_1^c, \|\cdot\|_\infty\big) \quad \lesssim \quad \frac{1}{c}\big(\log\frac{1}{\varepsilon}\big)^2.$$

*Proof.* We construct an $\varepsilon$-net of piecewise polynomials over $\mathbb{H}_1^c$.

Because all moments of the spectral measure $\mu_c$ are finite by (2.2), we can for any $h \in L_2(\mu_c)$ differentiate the function $\mathcal{F}_c h$ under the integral sign to find that $(\mathcal{F}_c h)^{(k)}(t) = \int(-i\lambda)^k e^{-it\lambda}h(\lambda) d\mu_c(\lambda)$. Consequently,

$$\left\|(\mathcal{F}_c h)^{(k)}\right\|_\infty^2 \leq \int |\lambda|^k |h(\lambda)| d\mu_c(\lambda)^2$$

$$\leq \int |h(\lambda)|^2 d\mu_c(\lambda) \int |\lambda|^{2k} d\mu_c(\lambda) = \|\mathcal{F}_c h\|_{\mathbb{H}^c}^2 \frac{\alpha_{2k}}{c^{2k}},$$

where $\alpha_k$ are the absolute moments of the spectral measure $\mu$. By Taylor's formula it follows that, for every $s, t \in [0,1]$,

$$(\mathcal{F}_c h)(t+s) = \sum_{j=0}^{k-1}(\mathcal{F}_c h)^{(j)}(t)\frac{s^j}{j!} + (R_k h)(t,s),$$

with remainder satisfying

$$\big|(R_k h)(t,s)\big| \leq \frac{1}{k!}\big\|(\mathcal{F}_c h)^{(k)}\big\|_\infty |s|^k \leq \frac{\sqrt{\alpha_{2k}}}{k!}\frac{|s|^k}{c^k},$$

for $\mathcal{F}_c h \in \mathbb{H}_1^c$.

For given $\varepsilon, d > 0$ choose $k \in \mathbb{N}$ such $\sqrt{\alpha_{2k}}(d/c)^k/k! \leq \varepsilon$. Consider the set $\mathcal{H}$ of functions

$$h(t) = \sum_{i=1}^{\lceil 1/d \rceil} 1_{((i-1)d), id]}(t) \sum_{j=0}^{k-1} \gamma_{i,j} \frac{(t-id)^j}{j!}, \qquad (2.3)$$

where $\gamma_{i,j}$ ranges over the grid $\{0, \pm\eta_j, \pm2\eta_j, \ldots\}$ intersected with the interval $\left[-\sqrt{\alpha_{2j}}/c^j, \sqrt{\alpha_{2j}}/c^j\right]$, for $\eta_j = \varepsilon\, j!/(d^j k)$. For every function $\mathcal{F}_c h \in \mathbb{H}_1^c$ and $i$ there exist points $\gamma_{i,j}$ in the grid such that

$$\sup_{(i-1)d < t \leq id} \Big| \sum_{j=0}^{k-1} (\mathcal{F}_c h)^{(j)}(id) \frac{(t-id)^j}{j!} - \sum_{j=0}^{k-1} \gamma_{i,j} \frac{(t-id)^j}{j!} \Big| \leq \sum_{j=0}^{k-1} \eta_j \frac{d^j}{j!} = \varepsilon.$$

The function $\sum_{j=0}^{k-1} (\mathcal{F}_c h)^{(j)}(id)(t-id)^j/j!$ is within uniform distance $\varepsilon$ of the function $\mathcal{F}_c h$ on the interval $((i-1)d, id]$. The preceding being true for every $i$ implies that the set $\mathcal{H}$ of piecewise polynomials (2.3) forms a $2\varepsilon$-net over $\mathbb{H}_1^c$ for the uniform norm on $(0,1]$, and hence the covering number $N(2\varepsilon, \mathbb{H}_1^c, \|\cdot\|_\infty)$ is bounded by the number of points in $\mathcal{H}$, which is equal to the number of different matrices $(\gamma_{i,j})$. The logarithm of this number can be bounded as

$$\log \#\mathcal{H} \leq \log \prod_{i=1}^{\lceil 1/d \rceil} \prod_{j=0}^{k-1} \Big(2\frac{\sqrt{\alpha_{2j}}/c^j}{\eta_j} + 1\Big) \leq \lceil \tfrac{1}{d} \rceil \sum_{j=0}^{k-1} \log\Big(2\frac{\sqrt{\alpha_{2j}}}{j!} \frac{k}{\varepsilon} \frac{d^j}{c^j} + 1\Big).$$

For any $x \geq 0$ and $j \geq 0$ we have $x^j \leq e^x \Gamma(j+1)$. Indeed,

$$e^x \Gamma(j+1) = \int_0^\infty e^{-(s-x)} s^j \, ds = \int_{-x}^\infty e^{-s}(s+x)^j \, ds \geq \int_0^\infty e^{-s} x^j \, ds.$$

Therefore, for any $\lambda$ we have $|\lambda|^{2j} \leq \delta^{-2j} \Gamma(2j+1) e^{\delta|\lambda|}$ and, consequently, with $K = \int e^{\delta|\lambda|} \, d\mu(\lambda)$,

$$\frac{\sqrt{\alpha_{2j}}}{j!} \Big(\frac{d}{c}\Big)^j \leq \frac{\sqrt{K\Gamma(2j+1)}}{\delta^j j!} \Big(\frac{d}{c}\Big)^j.$$

In view of Stirling's approximation $\Gamma(n+1) \asymp n^{n+1/2} e^{-n}$, the right-hand side is, up to a constant, for $j \geq 1$, equivalent to

$$\frac{\sqrt{K}(2j)^{j+1/4} e^{-j}\big(1 + o(1)\big)}{j^{j+1/2} e^{-j}\big(1+o(1)\big)} \Big(\frac{d}{\delta c}\Big)^j \lesssim \Big(\frac{2d}{\delta c}\Big)^j.$$

We choose $d < \delta c/2$ and $k \sim \log(1/\varepsilon)$ to reduce this expression for $j = k$ to a number smaller than $\varepsilon$. We have that $(d/c)^j \sqrt{\alpha_{2j}}/j!$ is bounded above uniformly in $j = 1, \ldots, k-1$, and hence

$$\log \#\mathcal{H} \lesssim \lceil \tfrac{1}{d} \rceil k \log \frac{k}{\varepsilon}.$$

With the indicated choices of $k$, this yields the bound given in the statement of the lemma. $\qquad\square$

If the spectral measure satisfies the stronger tail condition $\int \exp(\delta|\lambda|^r)\,\mu(d\lambda) < \infty$ for some $\delta > 0$ and $r > 1$, the elements of the RKHS can be extended to the entire complex plane and satisfy an exponential type restriction. In that case the results of Section 7.4 of [7] apply and can be used to improve the power 2 appearing in the entropy bound given by the lemma. In the statistical results, this improves the power of the logarithmic factors that we have in the rate of contraction results.

The following is a consequence of the preceding lemma and the well-known connection between the entropy of the unit ball of the RKHS and small ball probabilities, cf. Kuelbs and Li [8], Li and Linde [14].

**Theorem 2.4.** *Suppose the spectral measure satisfies (2.2) and $c \leq 1$. Then there exists an $\varepsilon_0 > 0$, independent of $c$, such that the rescaled process $W^c$ satisfies*

$$-\log \Pr\Big(\sup_{0 \leq t \leq 1} |W_t^c| \leq 2\varepsilon\Big) \lesssim \frac{1}{c}\Big(\log \frac{1}{c\varepsilon^2}\Big)^2$$

*for $\varepsilon \in (0, \varepsilon_0)$.*

*Proof.* Let $\psi_c(\varepsilon) = -\log \Pr\big(\sup_{0 \leq t \leq 1} |W_t^c| \leq \varepsilon\big)$ be the quantity of interest. The preceding lemma and Theorem 1.2 of [14] imply we have the crude bound

$$\psi_c(\varepsilon) \lesssim c^{-2/(2-\alpha)} \varepsilon^{-2\alpha/(2-\alpha)} \tag{2.4}$$

for every $\alpha \in (0, 2)$. According to the proof of Theorem 1.2 and the related Proposition 3.1 of [14], this bound holds for all $\varepsilon > 0$ satisfying

$$c \lesssim (\psi_c(\varepsilon/2))^{\alpha/2} \varepsilon^{-\alpha}.$$

We have $c \leq 1$ by assumption and hence

$$\psi_c(\varepsilon/2) = -\log \mathbb{P}\big(\sup_{0 \leq t \leq 1} |W_{t/c}| \leq \varepsilon/2\big) \geq -\log \mathbb{P}\big(\sup_{0 \leq t \leq 1} |W_t| \leq \varepsilon/2\big).$$

Since the right-hand side is independent of $c$, it follows that (2.4) holds for all $\varepsilon$ in an interval independent of $c$. The preceding lemma and Theorem 2 of [8] imply that for $\varepsilon$ small enough

$$\psi_c(2\varepsilon) \lesssim \frac{1}{c}\Big(\log \frac{\sqrt{\psi_c(\varepsilon)}}{\varepsilon}\Big)^2.$$

Again, inspection of the proof of the cited result of [8] shows that under our assumption $c \leq 1$, this holds for all $\varepsilon > 0$ in an interval independent of $c$. Combination of the preceding display with (2.4) now yields the statement of the theorem. $\square$

### 2.2. Multiply integrated Brownian motion

Consider a mean-zero Gaussian process $(W_t : t \geq 0)$ that is self-similar of order $\alpha$: the processes $(c^\alpha W_{t/c} : t \geq 0)$ and $(W_t : t \geq 0)$ are equal in distribution for every $c > 0$. The rescaled process $(W_{t/c} : 0 \leq t \leq 1)$ given in (1.1) is then equal in distribution to the process $(c^{-\alpha} W_t : 0 \leq t \leq 1)$, which means that for use as a prior distribution the rescaling of the time-axis is equivalent to a rescaling of the vertical axis. The rescaling has a simple effect on the reproducing kernel Hilbert space and small ball probability, but it has an interesting consequence.

We assume that the restriction $(W_t : 0 \leq t \leq 1)$ of the process to the unit interval and the rescaled process $W^c = (W_{t/c} : 0 \leq t \leq 1)$ can be viewed as Borel measurable maps in a separable Banach space $(\mathbb{B}, \|\cdot\|)$, and that the self-similarity can be understood in the sense that the Borel laws of these two processes are identical. The RKHS and small ball probability of the process $W^c = (W_{t/c} : 0 \leq t \leq 1)$ are then equal to these objects for the process $(c^{-\alpha} W_t : 0 \leq t \leq 1)$. Let $\mathbb{H}_W$ be the RHKS of the process $W$ (restricted to the unit interval) and let $\varphi_0(\varepsilon; W) = -\log \Pr(\|W\| \leq \varepsilon)$ be the exponent in its centered small ball probability. The following lemma is clear from the preceding.

**Lemma 2.5.** *The RKHS of the process $W^c$ is the set of functions $\mathbb{H}_W$ equipped with the norm $\|h\|_{W^c} = c^\alpha \|h\|_W$. The centered small ball exponent of $W^c$ satisfies $-\log \Pr(\|W^c\| \leq \varepsilon) = \varphi_0(c^\alpha \varepsilon; W)$.*

As an example consider the $k$-fold integrated Brownian motion. Define $I_{0+}^1 f$ as the function $t \mapsto \int_0^t f(s)\,ds$ and set $I_{0+}^k f = I_{0+}^1(I_{0+}^{k-1} f)$. Because Brownian motion $B$ is self-similar of index $1/2$, the process $W = I_{0+}^k B$ is self-similar of order $\alpha = k + 1/2$. We consider the restriction of this process to $[0,1]$ as a map in $C[0,1]$.

The fact that the integrated Brownian motion has $k$ derivatives at 0 equal to zero causes that the functions in its reproducing kernel Hilbert space satisfy similar constraints at 0. A better prior is obtained by adding an independent polynomial to the process. We consider the modified process

$$V_t^{c,a} = (I_{0+}^k B)_{t/c} + \frac{1}{\sqrt{a}} \sum_{i=0}^k Z_i \frac{t^i}{i!}, \qquad (2.5)$$

for scaling factors $c, a > 0$, $B$ a standard Brownian motion and independent standard normal variables $Z_0, \ldots, Z_k$, independent of $B$.

The following theorem gives a centered small deviation bound for the process $V^{c,a}$, and describes the approximation of smooth functions by elements of its RKHS $\mathbb{H}^{c,a}$.

**Theorem 2.6.** *Consider the process $V^{c,a}$ given in (2.5) as a map in $C[0,1]$. This process satisfies, for $\varepsilon > 0$ small enough,*

$$-\log \Pr\Big( \sup_{0 \leq t \leq 1} |V_t^{c,a}| \leq 2\varepsilon \Big) \lesssim \Big( \frac{1}{c^{k+1/2} \varepsilon} \Big)^{1/(k+1/2)} + k \log \frac{1}{\sqrt{a}\varepsilon}.$$

*Moreover, for $w \in C^\beta[0,1]$ and $\beta \leq k + 1$,*

$$\inf\big\{ \|h\|_{\mathbb{H}^{c,a}}^2 : \|h - w\|_\infty \leq \varepsilon \big\} \lesssim c^{2k+1} \Big( \frac{1}{\varepsilon} \Big)^{(2k+2-2\beta)/\beta} + a\Big( \frac{1}{\varepsilon} \Big)^{((2k-2\beta)/\beta)\vee 0}.$$

*Proof.* The assertion on the small ball probability follows easily from Theorem 2.1 of Li and Linde [13] on the small ball probability of integrated Brownian motion, and the fact that the added polynomial is independent and finite-dimensional.

By general arguments (e.g. [21], Section 10) we have that the reproducing kernel Hilbert space of the process (2.5) viewed as a map in $C[0, 1]$ is the Sobolev space $H^{k+1}[0, 1]$ of functions $h : [0, 1] \to \mathbb{R}$ that are $k$ times continuously differentiable with absolutely continuous $k$th derivative that is the integral of a function $h^{(k+1)} \in L_2[0, 1]$, equipped with the norm with square

$$\|h\|_{\mathbb{H}^{c,a}}^2 = c^{2k+1} \|h^{(k+1)}\|_2^2 + a \sum_{i=0}^{k} h^{(i)}(0)^2.$$

For a smooth function $\varphi$ and $\varphi_\sigma(x) = \varphi(x/\sigma)/\sigma$ its scaled version, the convolution $w * \varphi_\sigma$ is contained in the RKHS, with square norm

$$c^{2k+1} \int (w * \varphi_\sigma^{(k+1)})(x)^2 \, dx + a \sum_{i=0}^{k} (w * \varphi_\sigma)^{(i)}(0)^2.$$

If $\varphi$ is chosen such that $\int \varphi(x) \, dx = 1$ and with zero moments of orders $1, \ldots, k$, then the distance $\|w * \varphi_\sigma - w\|_\infty$ in the uniform norm can be seen to be of order $\sigma^\beta$. We choose $\sigma = \varepsilon^{1/\beta}$ and next evaluate the preceding display to be of the order as given in the theorem. $\qquad\square$

## 3. Posterior contraction rates

In this section we present the main results on posterior convergence rates using rescaled Gaussian process priors. We denote the posterior distribution based on a prior $\Pi_n$ and observations $X^{(n)}$ by $B \mapsto \Pi_n(B \mid X^{(n)})$.

We consider three different statistical settings: i.i.d. density estimation, classification, and fixed design regression. For any of these, the general theory developed in [20] gives results expressing posterior contraction rates in terms of the small deviations behaviour and the RKHS structure of the Gaussian prior. The results in the present section are obtained by combining these general results with the material of the preceding section.

Below we give complete proofs for the density estimation case. Since the other two cases are completely analogous we only explain the results briefly.

### 3.1. Density estimation

Suppose that we observe a random sample $X_1, \ldots, X_n$ from a positive density $\bar{p}_0$ on $[0, 1]$. A prior distribution on the set of positive densities can be defined structurally as $p_W$, for a Gaussian process $W = (W_t : t \in [0, 1])$ and $p_w$ the function defined by

$$p_w(t) = \frac{e^{w_t}}{\int_0^1 e^{w_t} \, dt}. \tag{3.1}$$

In the next two theorems we assume that the true density is $\alpha$-smooth in the sense that $\log \bar{p}_0 \in C^\alpha[0, 1]$. We show that if in this case we take for $W$ a suitably rescaled Gaussian process, we obtain a posterior that (perhaps up to logarithmic factors) contracts around the true density at the optimal minimax rate $n^{-\alpha/(1+2\alpha)}$.

The first result deals with rescaled smooth stationary processes.

**Theorem 3.1.** *Let $\alpha > 0$ be fixed. Let $W = (W_t : t \geq 0)$ be a centered, stationary Gaussian process with spectral measure $\mu$ satisfying condition (2.2) for some $\delta > 0$, and possessing a positive Lebesgue density. Let $W^n = (W_{t/c_n} : t \in [0, 1])$ be the rescaled version of $W$, for scaling constants $c_n \to 0$. Define the prior $\Pi_n$ structurally as $p_{W^n}$, with $p_w$ as in (3.1). Then if $\log \bar{p}_0 \in C^\alpha[0, 1]$, we have*

$$\mathbb{E}_0 \Pi_n(p : h(p, \bar{p}_0) > M\varepsilon_n \mid X_1, \ldots, X_n) \to 0$$

*for all $M$ large enough, where $\varepsilon_n = c_n^\alpha \vee (\log n)/\sqrt{nc_n}$ and $h$ is the Hellinger distance on densities. For*

$$c_n = \left( \frac{\log^2 n}{n} \right)^{\frac{1}{2\alpha+1}}$$

*this gives the rate $\varepsilon_n = (n/\log^2 n)^{-\frac{\alpha}{1+2\alpha}}$.*

*Proof.* By Theorem 3.1 of [20] we get the conclusion of the theorem as soon as we show that $\varphi_n(\varepsilon_n) \lesssim n\varepsilon_n^2$, where

$$\varphi_n(\varepsilon_n) = \inf_{h \in \mathbb{H}_n : \|h - \log \bar{p}_0\|_\infty < \varepsilon_n} \|h\|_{\mathbb{H}_n}^2 - \log \Pr(\|W^n\| < \varepsilon_n),$$

with $\mathbb{H}_n$ the RKHS of the rescaled process $W^n$. Hence, by Lemma 2.2 and Theorem 2.4, it suffices to verify that

$$\frac{1}{c_n} \left( \log \frac{1}{c_n \varepsilon_n^2} \right)^2 \lesssim n\varepsilon_n^2, \quad c_n^\alpha \leq \varepsilon_n \quad \text{and} \quad \frac{1}{c_n} \lesssim n\varepsilon_n^2.$$

It is easy to check that these relations indeed hold for $c_n$ and $\varepsilon_n$ as in the statement of the theorem. $\square$

The following theorem gives the analogous result for rescaled integrated Brownian motions.

**Theorem 3.2.** *For $\alpha > 0$ and $k \in \mathbb{N}_0$, let $V^n$ be the modified $k$-fold integrated Brownian motion defined in (2.5), with the scaling constant $c$ replaced by*

$$c_n = n^{\frac{\alpha - (k+1/2)}{(k+1/2)(1+2\alpha)}}$$

*and $a$ replaced by a sequence $a_n$ satisfying $a_n \leq n^{\frac{1+2\alpha-2k}{1+2\alpha}}$. Define the prior $\Pi_n$ structurally as $p_{V^n}$, with $p_w$ as in (3.1). Then if $\log \bar{p}_0 \in C^\alpha[0, 1]$ and $\alpha \leq k+1$, we have*

$$\mathbb{E}_0 \Pi_n(p : h(p, \bar{p}_0) > M\varepsilon_n \mid X_1, \ldots, X_n) \to 0$$

*for all $M$ large enough, where $\varepsilon_n = n^{-\frac{\alpha}{1+2\alpha}}$ and $h$ is the Hellinger distance on densities.*

*Proof.* By Theorem 3.1 of [20] and Theorem 2.6 it suffices to verify

$$\left(\frac{1}{c_n^{k+1/2}\varepsilon_n}\right)^{1/(k+1/2)} + k\log\frac{1}{\sqrt{a}\varepsilon_n} \lesssim n\varepsilon_n^2$$

and

$$c_n^{2k+1}\left(\frac{1}{\varepsilon_n}\right)^{(2k+2-2\beta)/\beta} + a\left(\frac{1}{\varepsilon_n}\right)^{(2k-2\beta)/\beta} \lesssim n\varepsilon_n^2.$$

For $c_n$, $\varepsilon_n$ and $a_n$ as in the theorem the left-hand sides of the displays are dominated by the first terms. Hence, it remains to check that

$$\left(\frac{1}{c_n^{k+1/2}\varepsilon_n}\right)^{1/(k+1/2)} \lesssim n\varepsilon_n^2$$

and

$$c_n^{2k+1}\left(\frac{1}{\varepsilon_n}\right)^{(2k+2-2\alpha)/\alpha} \lesssim n\varepsilon_n^2,$$

which is straightforward. ☐

### *3.2. Fixed design regression*

Suppose that we observe independent variables $Y_1, \ldots, Y_n$ following the regression model $Y_i = w_0(t_i) + e_i$ for unobservable $N(0, \sigma_0^2)$-distributed errors $e_i$ and fixed, known elements $t_1, \ldots, t_n$ of the unit interval. Consider estimating the regression function $w$.

As a prior on $w$ we use the Gaussian processes $W^n$ from Theorem 3.1 or $V^n$ from Theorem 3.2. If the standard deviation $\sigma_0$ of $e$ is not known, we also put a prior on $\sigma_0$, which we assume to be supported on a given interval $[a, b] \subset (0, \infty)$ with a Lebesgue density that is bounded away from zero.

A combination of Theorem 3.3 of [20] and the results of Section 2 then shows that if $w_0 \in C^\alpha[0, 1]$ and $\sigma_0 \in [a, b]$ the analogues of the statements of the Theorems for the density estimation case are true in this setting as well. We get the same rates of posterior contraction, and the statement of the theorems has to be replaced by

$$\mathbb{E}_0\Pi_n((w, \sigma) : \|w - w_0\|_n + |\sigma - \sigma_0| > M\varepsilon_n \,|\, Y_1, \ldots, Y_n) \to 0$$

for all $M$ large enough, where $\|f\|_n^2 = n^{-1}\sum f^2(t_i)$.

### *3.3. Classification*

Suppose that we observe a random sample of vectors $(X_1, Y_1), \ldots, (X_n, Y_n)$ from the distribution of $(X, Y)$, where $Y$ takes its values in the set $\{0, 1\}$ and $X$ takes its values in the unit interval. Consider estimating the binary regression function $\bar{f}_0(t) = \Pr(Y = 1 \,|\, X = t)$.

We construct a prior on the set of regression functions as $f_{W^n}$ (or $f_{V^n}$) for $W^n$ (or $V^n$) the Gaussian process from Theorem 3.1 (or 3.2) and $f_w$ the function $f_w(t) = \Psi(w_t)$, where $\Psi : \mathbb{R} \to (0, 1)$ is (for instance) the logistic distribution function.

Theorem 3.2 of [20] and the results of Section 2 imply that if $\Psi^{-1}(\bar{f}_0) \in C^\alpha[0, 1]$, the analogues of the statements of Theorems 3.1 and 3.2 hold in this setting. We get the same rates of posterior contraction in this case as well, the statement of the theorems has to be replaced by

$$\mathbb{E}_0 \Pi_n(f : \|f - \bar{f}_0\|_{G,2} > M\varepsilon_n \mid (X_1, Y_1), \ldots, (X_n, Y_n)) \to 0$$

for any sufficiently large constant $M$, where $\|f\|_{G,2}^2 = \int f^2(t)\,dG(t)$ and $G$ is the marginal distribution of $X$.